\documentclass[11pt]{amsart}

  \usepackage{epsf}

  \title{  Higher Lawrence configurations}

  \author{     Francisco Santos}
  \address{    Francisco Santos,
               Departamento de Matem\'aticas, Estad\'{\i}stica y
                   Computa\-ci\'on\\
               Universidad de Cantabria, E-39005, Santander, SPAIN.}
  \email{santos@matesco.unican.es}
  \author{     Bernd Sturmfels}
  \address{    Bernd Sturmfels,
               Department of Mathematics, University of California,
               Berkeley, CA 94720, USA.}
  \email{bernd@math.berkeley.edu}

  \thanks{The first author was supported by grant BFM2001--1153 of the  
Spanish
    Direcci\'on General de Investigaci\'on. The second author was
   partially supported by grants DMS-0200729 and
DMS-0138323 of the U.S.~National Science Foundation.}


  \newtheorem{theorem}{Theorem}
  \newtheorem{lemma}[theorem]{Lemma}

  \newtheorem{remark}[theorem]{Remark}
  \newtheorem{corollary}[theorem]{Corollary}

  \theoremstyle{definition}

  \newtheorem{example}[theorem]{Example}

  \newcommand{\naturals}{{\mathbb N}}
  \newcommand{\integers}{{\mathbb Z}}
  \newcommand{\reals}{{\mathbb R}}

  \newcommand{\e}{{\epsilon}}
  
  \newcommand{\C}{{\mathcal C}}
  
  \renewcommand{\L}{{\mathcal L}}
  \renewcommand{\P}{{\mathcal P}}
  \newcommand{\A}{{\mathcal A}}
  \newcommand{\B}{{\mathcal B}}

\begin{document}

  \begin{abstract}
Any configuration of lattice vectors gives rise to
a hierarchy of higher-dimensional configurations
which generalize the  Lawrence construction
in geometric combinatorics. We prove finiteness
results for the Markov bases, Graver bases and
face posets of these configurations, and we discuss
applications to the statistical theory of log-linear models.

  \end{abstract}

  \maketitle

  \section{\bf Introduction}
Fix a configuration $\A=\{a_1,\dots,a_n\}$ of lattice vectors  
spanning $\integers^d$,  and let $\,\L(\A)\subset \integers^n\,$ be
the lattice of linear relations on
$\,\A$. We introduce a
  hierarchy of configurations $\A^{(2)}, \A^{(3)}, \A^{(4)}, \dots$, as  
follows.
  The configuration
  $\A^{(r)}$ consists of $r\cdot n$ vectors in $\integers^{dr+n}=
  (\integers^{d}\otimes\integers^{r})\oplus\integers^{n}$, namely,
  \begin{equation}
  \A^{(r)}
\quad = \quad \bigl\{ \,
(a_i\otimes e_j) \oplus \e_i \,\,: \,i=1,\dots, n,\, j=1,\ldots,r  
\bigr\}
  \label{eqn:Ar}
  \end{equation}
  where $e_i$ and $e_j$ denote unit
vectors in $\integers^n$ and $\integers^r$ respectively.
The first object in this hierarchy is  $\A^{(2)}$,  which is
a configuration of $2n$ vectors isomorphic to the  \emph{Lawrence  
lifting}
$\Lambda(\A)$ of the given configuration  $\A$. See
\cite[\S 7]{Sturmfels-book} or \cite[\S 6.6]{Ziegler}.
  We call $\A^{(r)}$ the \emph{$r$-th Lawrence lifting of $\A$}.
In this paper we study the  Lawrence hierarchy $\,\A^{(2)}, \A^{(3)},
  \A^{(4)}, \dots \,$ from the perspective of toric algebra,
geometric combinatorics and applications to statistics \cite{Diaconis}.

The $r$-th Lawrence lifting  $\A^{(r)}$ is characterized as the
configuration whose linear relations are $r$-tuples of
linear relations on $\A$  that sum to zero. Indeed,
the lattice of linear relations on $\A^{(r)}$
has rank $(r-1)(n-d)$ and equals
$$
  \L(\A^{(r)}) \,\, = \,\,
  \left\{     (u^{(1)} \! , u^{(2)} \!, \dots, u^{(r)}) \in  
(\integers^n)^r
      :
      u^{(i)} \in \L(\A)\ \forall i,\  \sum u^{(i)}=0   \right\}.
$$
We think of the elements of $\L(\A^{(r)})$ as integer
  $r\times n$-tables whose column sums are zero and whose $\A$-weighted
row sums are zero. The \emph{type} of such a table
is the number of non-zero row  vectors $u^{(i)}$.
  Given any basis $\,\{b^{(1)}, \dots, b^{(n-d)}\,\}
\subset \integers^n\,$ for the lattice $\,\L(\A)\,$ of linear relations
on $\,\A$,  it is easy to derive a lattice basis of $\L(\A^{(r)})$
consisting of tables of type $2$. For instance, take tables with first  
row
some $\, b^{(i)}\,$ and some other row $\,-b^{(i)}$.

  In toric algebra and its statistics applications we are
interested in larger subsets of   $\L(\A^{(r)})$ which generate
the lattice in a stronger sense. A \emph{Markov basis} of $\A^{(r)}$
is a finite subset of $\L(\A^{(r)})$ which corresponds to a minimal
set of generators of the toric  ideal $I_{\A^{(r)}}$
as in \cite[\S 4]{Sturmfels-book}, or, equivalently,
to a minimal set of moves which connects any
two nonnegative integer $r\times n$-tables that have the same
column sums and the  same $\A$-weighted row sums  \cite{AS1},
\cite{Diaconis}, \cite{GMS}, \cite{HS}.
We prove that the Markov bases stabilize for $r \gg 0$.

  \begin{theorem}
  \label{first}
For any configuration $\A=\{a_1,\dots,a_n\}$ in $\integers^d$,
there exists a constant $m=m(\A)$ such that any
higher Lawrence lifting $\A^{(r)}$, for any $r \geq 2$, 
has a Markov basis
consisting of tables having type at most $m$.
  \end{theorem}

  We call the minimum value $m(\A)$ the \emph{Markov complexity} of $\A$.
This paper was inspired by recent work of
the statisticians Aoki and Takemura  \cite{AS1}. Their result
states, in our notation, that  the product of two triangles
  \[
  \Delta_2\times \Delta_2
\quad  = \quad
  \{\,e_i\oplus e_j\in \integers^3\oplus \integers^3 : 1\le i,j\le 3  
\,\}\,
  \]
has Markov complexity $5$.
Indeed, $\L((\Delta_2\times \Delta_2)^{(r)})$  consists of
integer $3\times 3\times r$-tables with zero line sums in the three  
directions.
These are all possible moves for the {\em no-three-way interaction  
model} \cite{Diaconis}. Our Theorem \ref{first} implies:

  \begin{corollary}
  \label{coro:threeway}
  For any positive integers $p$ and $q$
there exists an integer $m$ such that the
Markov basis for $p\times q\times r$-tables (in the no three-way  
interaction
model, for arbitrary $r$) consists of tables of format $p\times q\times  
m'$ with $m'\le m$. \end{corollary}

We often use the phrase ``the Markov basis''
instead of ``a Markov basis''. The definite
article  is
justified because the  minimal generating set 
of a homogeneous toric ideal is unique up
to minor combinatorial modifications.

We prove Theorem \ref{first} by providing
an explicit upper bound for $m(\A)$.
Recall that the  \emph{Graver basis} of
$\A$ is the set of minimal
elements in $\L(\A) \backslash \{0\}$,
where $\integers^n$ is partially ordered by
setting $a \leq b\,$ if  $\,b \,$ is the
\emph{conformal} sum of $a$ and $b-a$.
This condition means that,
for every $i \in \{1,\ldots,n\}$, either
$ 0 \leq a_i \leq b_i $ or $ 0 \geq a_i \geq b_i$ holds.
The Graver basis is unique, finite, and contains
Markov bases for all subconfigurations of $\A$.
See \cite{Sturmfels-book} for bounds, algorithms and many details.
We define the {\em Graver complexity} $g(\A)$ to be the maximum type
of any table in the Graver basis of some
higher Lawrence lifting $\A^{(r)}$. Clearly, $m(\A) \leq g(\A)$.
We now state our main result. The phrase
``the Graver basis of the Graver basis''
is not a typo but it is the punchline.
We regard the elements in the
 Graver basis of $\A$ as the columns
of some big matrix and then we
compute the Graver basis of that big matrix.

  \begin{theorem}
  \label{main}
The Graver complexity $g(\A)$ of a configuration $\A$
is the maximum $1$-norm
of any element in the Graver basis of the Graver basis of~$\A$.
  \end{theorem}

This paper is organized as follows.
In Section 2 we present a few examples
to illustrate the notions of Markov complexity
and Graver complexity. The proof of Theorem \ref{main}
(and hence of Theorem \ref{first}) will be given in Section 3.
Section 4 deals
with applications to statistics. We show that
if $\A$ is any log-linear hierarchical model
(in the notation of \cite{GMS}, \cite{HS})
then $\A^{(r)}$ is the corresponding
{\em logit model} (in the sense of \cite[\S VII]{Chr}, \cite[\S 6]{Fie})
where the {\em response variable}
has $r$ levels. Thus Theorem \ref{first} implies
the existence of a finite Markov basis for
logit models where the response variable
has an unspecified number of levels.
   In Section 5 we prove an analogue of
Theorem \ref{main} for circuits,
and we examine the convex polytopes arising from
higher Lawrence liftings $\A^{(r)}$.

\section{\bf Examples}
\label{sec:examples}

The first three examples below show that the Markov and Graver  
complexities of a
configuration may coincide or differ a lot. After this we work out in
detail the {\em twisted cubic curve}, a familiar
example in toric algebra.

\begin{example}
  \label{exm:lawrence}
  Let $\A$ be any configuration and $\Lambda(\A)=\A^{(2)}$ its usual
Lawrence  lifting. Then
$\,\Lambda(\A)^{(r)}$ is the Lawrence lifting  of $\A^{(r)}$,
  since $(\A^{(r)})^{(s)}= (\A^{(s)})^{(r)}$ for all $r,s$.
By~\cite[Theorem~7.1]{Sturmfels-book}, the Markov and Graver complexity
of $\Lambda(\A)$ coincide. They are equal to the Graver complexity
of $\A$.
  \end{example}

\begin{example}
\label{exm:productofsimplices}
Let $\A=\{1,\dots,1\}$ consist of $n$ copies of the vector $1$
in $\integers^1$. The $r$-th Lawrence lifting $\A^{(r)}$
is the product of two simplices $\Delta^{n-1}\times \Delta^{r-1}$.
In statistics, this corresponds to two-dimensional tables of size
$n \times r$.
The Graver basis of $\A^{(r)}$ consists of the
circuits in the complete bipartite graph $K_{n,r}$,
and the Markov basis consists of circuits
which fit in a subgraph $K_{n,2}$.
The Graver complexity of $\A$ is $n$,
and  the Markov complexity of $\A$ is $2$.
\end{example}

\begin{example}
\label{exm:pqr}
Take $d=1$, $n=3$ and  $\,\A=\{k,l,m\}$, where $k$, $l$ and $m$ are
pairwise relatively prime. Using Theorem \ref{main}, it can be shown
that the Graver complexity $g(\A)$ equals $k+l+m$. We invite
the reader to write down the Graver basis element of type $k+l+m$ for
$\A^{(k+l+m)}$. It would be interesting to find a formula,
in terms of $k,l$ and $m$, for the Markov complexity $m(\A)$.
\end{example}

\begin{example}
\label{exm:twistedcubic} {\sl (Twisted Cubic) }
Let $d = 2 $ and $n = 4$ and fix the configuration
\begin{equation}
\label{twiscubmatrix}
  \A  \quad = \quad
\begin{pmatrix}
3 & 2 & 1 & 0 \\
0 & 1 & 2 & 3 \\
\end{pmatrix}
\end{equation}
The corresponding statistical model is
{\em Poisson regression} with four levels.
The toric ideal $I_\A$ of this configuration
consists of the
algebraic relations among the four cubic
monomials in two unknowns $s$ and $t$:
\begin{eqnarray*}
& I_\A \quad = \quad
  \langle \,
x_1 x_3 - x_2^2, \,
x_1 x_4 - x_2 x_3, \,
x_2 x_3 - x_3^2 \, \rangle \qquad = \\
&  {\rm kernel}
\bigl( k[x_1,x_2,x_3,x_4] \rightarrow k[s,t], \,
x_1 \mapsto s^3, \,
x_2 \mapsto s^2 t, \,
x_3 \mapsto s t^2, \,
x_4 \mapsto t^3 \bigr) .
\end{eqnarray*}
The Markov basis of $\A$ is the set of three
vectors $(1,-2,1,0)$, $(1,-1,-1,1)$, and $(0,1,-2,1)$
corresponding to the minimal generators of $I_\A$.
The Graver basis of $\A$ has two additional elements,
namely $(1,0,-3,2)$ and $(2,-3,0,1)$.

The ``classical'' Lawrence lifting is isomorphic to the
eight column vectors of
$$ \A^{(2)} \quad = \quad
\begin{pmatrix}
\A & {\bf 0} \\
{\bf 0} & \A \\
{\bf 1} & {\bf 1}
\end{pmatrix},
$$
where ${\bf 1}$ is the identity matrix of size $4\times 4$.
This $8 \times 8$-matrix has rank $6$.
Its kernel $\, \L(\A^{(2)}) \,$ is a rank $2$ lattice
whose elements are identified with
  $2 \times 4$-integer tables $T$
with $\, (1 \,\,\, 1) \cdot T\, = \, {\bf 0} \,$
and $\,T \cdot \A^t \, = \, {\bf 0} $. It is spanned by
$$
  \begin{pmatrix} 1 & -2 & 1 & 0 \\ -1 & 2 & -1 & 0 \end{pmatrix}
\quad \hbox{and} \quad
  \begin{pmatrix} 1 & -1 & -1 & 1 \\ -1 & 1 & 1 & -1 \end{pmatrix}.
$$
By \cite[Theorem~7.1]{Sturmfels-book},
the Markov basis of $\A^{(2)}$ equals the
Graver basis of $\A^{(2)}$. It consists of
the five tables constructed from the Graver basis of $\A$:
  \begin{eqnarray*}
& \begin{pmatrix} 1 & -2 & 1 & 0 \\ -1 & 2 & -1 & 0 \end{pmatrix},
  \begin{pmatrix} 1 & -1 & -1 & 1 \\ -1 & 1 & 1 & -1 \end{pmatrix},
  \begin{pmatrix} 0 & 1 & -2 & 1 \\ 0 & -1 & 2 & -1 \end{pmatrix}, \\
& \begin{pmatrix} 1 & 0 & -3 & 2 \\ -1 & 0 & 3 & -2 \end{pmatrix},
  \begin{pmatrix} 2 & -3 & 0 & 1 \\ -2 & 3 & 0 & - 1 \end{pmatrix}.
  \end{eqnarray*}
The third Lawrence lifting consists of the columns of
the $10 \times 12$-matrix
\begin{equation}
\label{Athree}
\A^{(3)} \quad = \quad
\begin{pmatrix}
\A & {\bf 0} & {\bf 0} \\
{\bf 0} & \A  & {\bf 0}\\
{\bf 0} & {\bf 0} & \A \\
{\bf 1} & {\bf 1}& {\bf 1}
\end{pmatrix}.
\end{equation}
Its kernel $\, \L(\A^{(3)}) \,$ is the rank $4$ lattice
consisting of $3 \times 4$-integer tables $T$
with $\, (1 \,\,\,1 \,\,\, 1) \cdot T\, = \, {\bf 0} \,$
and $\,T \cdot \A^t \, = \, {\bf 0} $. The
Markov basis of $\A^{(3)}$ has $21$ tables. Fifteen
of them are gotten from the five tables above by adding a row of zeros.
The other six Markov  basis elements are  row permutations of
\begin{equation}
  \begin{pmatrix}
  0 & 1 & -2 & 1 \\
  1 & -2 & 1 & 0 \\
  -1 & 1 & 1 & -1 \\
  \end{pmatrix}.
\label{eqn:markovtype3}
\end{equation}
It can be checked that no new Markov basis elements are needed
for $\, \A^{(r)}$, $ r \geq 4 $. Any
two $r \times 4$-tables of non-negative integers
which have the same column sums and the same
$\A$-weighted row sums can be connected by
the known moves involving only two or three of the rows.
Equivalently:
\begin{remark}
The twisted cubic curve $\A$ has Markov complexity  $m(\A) = 3$.
\end{remark}

We next discuss the Graver complexity of $\A$.
The Graver basis of $\A^{(3)}$ consists of $87$ tables.
Every other table in $ \, \L(\A^{(r)}) \,$ can be expressed as
an $\naturals$-linear combination of these $87$
without cancellation in any coordinate.
In addition to the $21$ Markov basis  elements,
the Graver basis of $\A^{(3)}$ contains the
following $66$ tables which come in $6$ symmetry classes
(with respect to permutations of the three rows and mirror
reflection of the columns):
\begin{eqnarray*} &
  \hbox{Class 1 ($12$ tables, degree $7$)}:\qquad\qquad\qquad
  \begin{pmatrix}
  -1 & 1 & 1 & -1 \\  0 & -1 & 2 & -1 \\  1 & 0 & -3 & 2 \\
  \end{pmatrix}
\\ &   \hbox{Class 2 ($12$ tables, degree $9$, circuit)}:\qquad\qquad
  \begin{pmatrix}
  0 & -2 & 4 & -2 \\  -1 & 2 & -1 & 0\\  1 & 0 & -3 & 2 \\
  \end{pmatrix} \\
  &
\hbox{Class 3 ($12$ tables, degree $9$)}: \qquad\qquad
  \begin{pmatrix}
  -2 & 2 & 2 & -2 \\  1 & -2 & 1 & 0 \\  1 & 0 & -3 & 2 \\
  \end{pmatrix}
\end{eqnarray*}
\begin{eqnarray*}
  &
  \hbox{Class 4 ($6$ tables, degree $10$)}: \qquad\qquad
  \begin{pmatrix}
  -2 & 3 & 0 & -1 \\  1 & -3 & 3 & -1 \\  1 & 0 & -3 & 2 \\
  \end{pmatrix}
\\ &
  \hbox{Class 5 ($12$ tables, degree $12$)}:\qquad\qquad
  \begin{pmatrix}
  -3 & 4 & 1 & -2 \\  2 & -4 & 2 & 0 \\  1 & 0 & -3 & 2 \\
  \end{pmatrix}
\\ &
  \hbox{Class 6 ($12$ tables, degree $15$, circuit)}: \qquad\qquad
  \begin{pmatrix}
  -4 & 6 & 0 & -2 \\  3 & -6 & 3 & 0 \\  1 & 0 & -3 & 2 \\
  \end{pmatrix}
\end{eqnarray*}
Here
``degree'' refers to the total degree of the associated
binomial, and
``circuit'' means that the table has minimal support
with respect to inclusion \cite[\S 4]{Sturmfels-book}.
For instance, the binomial
of degree $15$ for the table in Class 6  is
$$
x_{12}^6 x_{21}^3 x_{23}^3 x_{31} x_{34}^2 \,\, - \,\,
x_{11}^4 x_{14}^2 x_{22}^6 x_{33}^3 $$
The Graver bases of $\A^{(4)}$ has $240$ elements of type four, and
hence it has
  \[
  240 \,+\, \binom{4}{3} \cdot 87  \, +\, \binom{4}{2}
\cdot 5  \quad = \quad 558
  \]
elements in total. We similarly compute the
Graver bases for the higher Lawrence liftings
$\A^{(5)}$, $ \A^{(6)}$,  $ \ldots $,
for instance, using Hemmecke's  program {\tt 4ti2} \cite{Hem}.
The Graver basis of $\A^{(6)}$ contains the following
table of type $6$:
  \begin{equation}
\label{sixbyfour}
  \hbox{($120$ tables, degree $15$, type 6) }\qquad\qquad
  \begin{pmatrix}
  -2 & 3 & 0 & -1 \\  -2 & 3 & 0 & -1 \\
  1 & -2 & 1 & 0 \\   1 & -2 & 1 & 0 \\   1 & -2 & 1 & 0 \\
  1 & 0 & -3 & 2 \\
  \end{pmatrix}
  \end{equation}
The Graver basis element (\ref{sixbyfour}) shows
that the Graver complexity of $\A$ is at least six.
Using Theorem \ref{main}, we can check that this is
the correct bound.

\begin{remark}
\label{rem:twistedGraver}
The twisted cubic curve $\A$ has Graver complexity
$g(\A) = 6$.
\end{remark}
\end{example}

  \section{\bf Proofs}
  \label{sec:proofs}

We first note that Theorem \ref{main} implies
Theorem \ref{first}, and hence also Corollary \ref{coro:threeway}.
The point is that the Graver basis of a toric ideal contains
a subset of minimal generators (i.e.~a Markov basis), and therefore
$m(\A) \leq g(\A)$. So, in order to show that
$m(\A)$ is finite, it suffices to show that
$g(\A)$ is finite.

To derive the exact formula for $g(\A)$  given
in Theorem  \ref{main}, we begin with the observation
that Graver basis elements of $\A$ are those
vectors $a$ in $\L(\A) \backslash \{0\}$ that cannot be decomposed
as a conformal sum $a=b+c$ with $b,c\in\L(\A)\backslash \{0\}$.
{\em Conformal} means $|b_i+c_i| = |b_i| + |c_i|$ for all $i$.

\begin{lemma}
\label{lem:decompose}
Let $u$ be a Graver basis element of $\A^{(r)}$ and
suppose that one of its rows, say $ a = u^{(i)}$,
has a conformal decomposition
$\,a =a_1+\cdots + a_k$, where the $a_i$'s are in $\L(\A)$.
Then the table $u'$,  gotten by removing the row $a$ from $u$ and
inserting the rows $a_1,\dots,a_k$,
is in the Graver basis of $\A^{(r+k-1)}$.
\end{lemma}

\begin{proof}
If $u'$ is not in the Graver basis, then it has a non-trivial
conformal decomposition
$u'=v' + w'$ with $v',w'\in\L(\A^{(r+k-1)})$.
Then, adding up the relevant $k$ rows of $v'$ to become a single row,
and the same for $w'$,
we get two tables $v,w\in\L(\A^{(r)})$ and a
non-trivial conformal decomposition
$u=v+w$ which proves that $u$ is not in the Graver basis
of $\A^{(r)}$ either.
\end{proof}

\begin{corollary}
\label{crucial}
Every Graver basis element $u$ of some $\A^{(r)}$ can be obtained by
conformal addition of rows from a Graver basis element $u'$
of some $\A^{(s)}$ which has the property
that each row of $u'$ lies in the Graver basis of $\A$.
\end{corollary}

Note that the implication of Lemma \ref{lem:decompose}
works only in one direction.
If $u$ is Graver then $u'$ is Graver,
but the converse is generally not true.

\begin{example}
\label{exm:repeatedrows}
Let $\A=\{1, 2, 1\}$. The first of the following
two tables is in the Graver basis of $\A^{(4)}$ but
the second is not in the Graver basis of $\A^{(3)}$.
\begin{equation}
u'\,=\,
\begin{pmatrix}
  0 & -1 &  2 \\
  2 & -1 &  0 \\
-1 &  1 & -1 \\
-1 &  1 & -1 \\
  \end{pmatrix},
\qquad
u\,=\,
\begin{pmatrix}
  0 & -1 &  2 \\
  2 & -1 &  0 \\
-2 &  2 & -2 \\
\end{pmatrix}
.
\label{eqn:repeatedrows}
\end{equation}
\end{example}

\vskip .1cm

\noindent {\sl Proof of Theorem \ref{main}: }
Let  $\B=\{b_1,\dots,b_k\}$ be the Graver basis of $\A$.
Corollary \ref{crucial} tells us that in computing
the Graver complexity  $g(\A)$ we only need to consider
tables all of whose rows lie in $\B$.
Let $\,u \,=\, (u^{(1)}, \ldots,u^{(r)}) \,\in \,\B^r \,$
be such a table, for $r \geq 3$, and suppose
that $u^{(i)} \not= - u^{(j)}$ for all $i,j$.
We define $\psi_u$ to be the integer vector of length $k$
whose $i$-th entry counts
(with sign) how many times $b_i$ appears as a row in $u$.
Then the $1$-norm of the vector $\psi_u$ equals
the number $r$, which is the type of the table $u$. Hence
the following claim will imply Theorem \ref{main}:
{\em The table $u$ is in the Graver basis of $\A^{(r)}$ if and
only if the vector $\psi_u$ is in the Graver basis of $\B$}.

To prove this claim, first suppose that
$\psi_u$ is not in the Graver basis of $\B$. Any
conformal decomposition of $\psi_u$  provides a
conformal decomposition of $u$ (into tables of smaller type),
so that $u$ is not in the Graver basis of $\A^{(r)}$. For the
converse, we note that any conformal decomposition of $u$
arises in this manner from some conformal decomposition
of $\psi_u$, because single rows of $u$ admit no conformal  
decomposition.
Hence any  non-trivial conformal decomposition of $u$
gives a  non-trivial conformal decomposition of $\psi_u$.
\qed

\vskip .2cm

It is instructive to examine the
proof of Theorem \ref{main} for each of
the examples discussed in Section 2.
For instance, if $\A$  is the twisted cubic
in  (\ref{twiscubmatrix}) then
$$
  \B  \quad = \quad
\begin{pmatrix}
  1 &  2 &  1 & 1 & 0 \\
-2 & -3 & -1 & 0 & 1 \\
  1 &  0 & -1 &-3 &-2 \\
  0 &  1 &  1 & 2 & 1 \\
\end{pmatrix}.
$$
The Graver basis of $\B$ consists of $13$ vectors,
ten of which are the circuits. The vector of
maximum $1$-norm among these $13$ vectors occurs for
$$ \psi_u \quad = \quad \bigl(\,3,\,-2,\,0,\,1,\,0\,\bigr), $$
the vector associated with  the $6 \times 4$-table $u$
in (\ref{sixbyfour}).

We can now derive a bound for $g(\A)$ in terms of $n$, $d$ and the  
maximum
size of the entries in $\A$, which we denote $s$.
Theorem 4.7 in \cite{Sturmfels-book} says that
the maximum $1$-norm of the vectors in $\B$ is at most
$(d+1)(n-d)D(\A)$, where $D(\A)\le (ds)^{d/2}$ is the maximum absolute  
value
among the full-dimensional minors of $\A$.
This implies bounds for the cardinality  $N$ of $\B$
and the maximum size $D(\B)$ of a subdeterminant of $\B$. For example:
\[
N\le (2(d+1)(n-d)D(\A) )^n,
\qquad
D(\B) \le \left((d+1)(n-d)D(\A)\right)^{n-d}
\]
Since $\B$ has dimension $n-d$, the same theorem cited above implies
\[
g(\A) \,\,\le\,\, (n-d+1)(N-(n-d)) D(\B) \,\,\le \,\,
n \left(2(d+1)(n-d)D(\A)\right)^{2n-d}.
\]

\begin{remark}
\rm
The finiteness of $g(\A)$ can also be
derived from a result about partially ordered
sets (posets) proved in 1952 by Higman \cite[Theorem 4.2]{Hig}.
We briefly present this approach which was suggested to us by
Matthias Aschenbrenner.
For any poset $S$, we can define a new poset
$\widetilde S$  as follows. The elements of $\widetilde S$
are the finite multisubsets of $S$, and the order is
  \begin{equation}
  \label{eqn:injectiveorder}
  V \le V'
  \quad \Leftrightarrow \quad
  \exists \,f: V\to V' \hbox{ injective and with }
\,\, \forall\, v \in V \,:\, v \le f(v)
  \end{equation}
A poset $S$ is said to be {\em Noetherian}
if every  non-empty subset of $S$  has
at least one, but at most finitely many minimal elements.
  Higman proved that  if $S$ is a Noetherian poset
then $\widetilde S$ is also Noetherian.
In his paper \cite{Hig}, he attributes this result to an earlier  
unpublished
manuscript of Erd\"os and Rado.

We apply this to the poset $\,S \,= \, \integers^n $,
defined as in the introduction:
$$ a \leq b \quad  \iff \quad
\hbox{for all $i \in \{1,\ldots,n\}$:
$ \,\,0 \leq a_i \leq b_i $ or $ 0 \geq a_i \geq b_i$}. $$
The poset $\,\integers^n \,$ is known to be Noetherian.
The  poset $\,\widetilde{\integers^n}\,$
consists of all finite multisubsets of $\integers^n$.
Higman's result implies that $\,\widetilde{\integers^n}\,$
is Noetherian.

There is a canonical map $\phi_r$ from the lattice $(\integers^n)^r$
of $r \times n$-tables to $\widetilde{\integers^n}$.
This map takes
$\,u\, = \, (u^{(1)} \! , u^{(2)} \!, \dots, u^{(r)}) \in  
(\integers^n)^r \,$
to the multiset of its non-zero row vectors
$\, \phi_r(u)\, = \,\bigl\{u^{(1)} \! ,
  u^{(2)} \!, \dots, u^{(r)}\bigr\}\backslash \{0\}$.
The union of the images of the maps $\phi_r$, as
$r$ ranges over $\naturals$, is the following subset of
$\widetilde{\integers^n}$:
$$ {\P (\A)}\,\, = \,\, \bigl\{ \,
V \in \widetilde{\integers^n} \,:\,\,
\hbox{the elements of $V$ lie in $\L(\A)\backslash\{0\}$ and
sum to zero} \,\bigr\}. $$
\begin{corollary}
The infinite set $\,{\P (\A)}\,$ has only finitely many
minimal elements, with the partial order induced
from $\widetilde{\integers^n}$.
The Graver complexity $g(\A)$ is the maximum
of their cardinalities.
\end{corollary}

\begin{proof}
The first assertion follows from Higman's result, which implies that
$\widetilde{\integers^n}$ is Noetherian. For the last assertion,
just observe that
an $r \times n$-table $u \in \L(\A^{(r)})$ lies in the Graver basis
of $\A^{(r)}$ if and only if the multiset
$\phi_r(u)$ is minimal in
$\,{\P (\A)} $, and
the type of $u$ is the cardinality of
$\phi_r(u)$.
\end{proof}

%

\end{remark}

  \section{\bf Statistics}

In this section we apply our results on higher Lawrence configurations
to the statistical context of log-linear models. We consider  
hierarchical
log-linear models for $m$-dimensional contingency tables. Such a model
is specified by a collection $\Delta$ of subsets of $\{1,2,\ldots,m\}$.
The standard notation for $\Delta$, used in the books of
Christensen \cite{Chr},  Fienberg \cite{Fie}
and other texts on cross-classified data,
is a string of brackets each containing the elements of a subset in  
$\Delta$.
For instance, the {\em four-cycle model} for $4$-dimensional tables
is $ \, \Delta = \bigl\{ \{1,2\} , \{2,3\}, \{3,4\}, \{4,1\} \bigr\}$,  
or,
in standard notation, $\, \Delta = [1 2][2 3][3 4][41]$.
(This is the smallest {\em non-decomposable graphical model}.)

If the format of the table is specified, say $\, r_1 \times r_2 \times
\cdots \times r_m$, then the model $\Delta$ is represented in toric
algebra by a configuration $\A$ as above, where $\, n = r_1 r_2 \cdots  
r_m \,$
and $d$ is the sum of the products $\, r_{i_1} r_{i_2} \cdots r_{i_p}  
\,$
where $\{ i_1,i_2, \ldots,i_p \}$ runs over $\Delta$.
  Each coordinate in $\A$ is either $0$ or $1$. For instance,
the four-cycle model $\, \Delta = [1 2][2 3][3 4][4 1]\,$ for
$2 \times 2 \times 2 \times 2$-tables is represented by
(the columns of) a $16 \times 16$-matrix with $0-1$-entries.
Here the lattice $\L ( \A )$ consists of all $2 \times 2 \times
2 \times 2$-integer tables whose
$[12]$-margins, $[23]$-margins, $[34]$-margins
and $[41]$-margins are zero.
See~\cite{GMS}, \cite{HS}.

The passage from the matrix $\A$ to its Lawrence lifting
$\A^{(r)}$ has the following statistical interpretation. Think of the
$m$ given random variables as {\em explanatory variables}, and consider
an additional $(m \! + \! 1)^{\rm st}$ random variable, the
{\em response variable}, which has $r$ levels. From the
model $\Delta$ for $m$-dimensional tables, we
construct the following model for $(m \! + \! 1)$-dimensional tables:
$$ \Delta_{\rm logit} \quad = \quad
\bigl\{ \, \{1,2,\ldots,m\} \, \bigr\} \,\, \cup \,\,
\bigl\{ \sigma \cup \{m+1\} \,\, :\,\, \sigma \in \Delta \,\bigr\}. $$
This is the {\em logit model} described in \cite[\S VII.1]{Chr}.
For example, if $\Delta$ is the four-cycle model and
the index ``$5$'' indicates
the additional response variable then, in standard notation,
$\, \Delta_{\rm logit} \,  =\,  [125][235][345][415][1234]$.
We shall prove that the passage from a log-linear model $\Delta$
to the associated logit model $\Delta_{\rm logit}$ is described
in toric algebra precisely by the Lawrence hierarchy.

\begin{theorem}
If $\A$ represents a hierarchical log-linear model $\Delta$
for $r_1 \times \cdots \times r_m$-tables then $\A^{(r)}$
represents the logit model $\Delta_{\rm logit}$ for
$r_1 \times \cdots \times r_m \times r$-tables.
\end{theorem}

\begin{proof}
We think of an $r_1 \times \cdots \times r_m \times r$-table
as a two-dimensional matrix with
$r_1 r_2 \cdots r_m$ columns and $r$ rows.
Computing the $\A$-weighted row sums of such a matrix
means computing the   $\,(\sigma \cup \{m+1\})$-marginals
for any $\sigma \in \Delta$. Computing the column sums
of such a matrix means computing the $\{1,2,\ldots,m\}$-marginals
of the $r_1 \times \cdots \times r_m \times r$-table. Thus
$\, \L (\A^{(r)}) \,$ is identified with the lattice of integer
$r_1 \times \cdots \times r_m \times r$-tables whose margins
in the model $\Delta_{\rm logit}$ are zero. This is precisely the claim.
\end{proof}

   From Theorems \ref{first} and \ref{main},
we obtain the following corollary.

\smallskip

\begin{corollary}
\label{struc}
Consider a logit model $\Delta_{\rm logit}$ where the
numbers $r_1,\ldots,r_m$ of levels of the explanatory
variables are fixed, and the number $r$
of levels of the response variable is allowed to increase.
Then there exists a finite Markov basis which is
independent of $r$, and independent of possible  structural zeros.
\end{corollary}

We need to explain the last subclause.
Imposing {\em structural zeros} in
the model $\A^{(r)}$ means to consider
the toric model defined by a subconfiguration
$\C \subset \A^{(r)}$. The Graver basis
of $\A^{(r)}$ is a universal Gr\"obner basis
\cite[\S 7]{Sturmfels-book}, and hence it
contains generators for all  elimination ideals.
This implies:

\begin{remark}
The Graver basis of $\A^{(r)}$ contains
a Markov basis for any  subconfiguration $\C
\subset \A^{(r)}$. (It works
even if structural zeros are imposed.)
\end{remark}

Hence to get the last assertion in Corollary \ref{struc},
one takes the Graver basis of $\A^{(r)}$ for $r \gg 0$.
The prototype of such a finiteness result   was obtained by
Aoki and Takamura in \cite{AS1}.
They considered the no-three-way interaction
model for three-dimensional contingency tables. This is the
logit model
$$ \Delta_{\rm logit} \quad = \quad [12][13][23] $$
derived from the most classical independence model
$$ \Delta \quad = \quad [1][2] $$
for two-dimensional tables.
If $r_1 = 2 \leq r_2$ then it was known from \cite{Diaconis}
that the Markov basis stabilizes for $r \geq r_2$. Aoki
and Takamura \cite{AS1} considered the case $r_1 = r_2 = 3$,
and they constructed the Markov basis which stabilizes
for $r \geq 5$. Using Theorem \ref{main} and
  Hemmecke's  program {\tt 4ti2} \cite{Hem},
  we found that the Graver basis
for $3 \times 3 \times r$-tables stabilizes
for $r \geq 9$. In symbols,
$$ g(\Delta_2 \times \Delta_2) \quad = \quad 9 .$$
An element of the Graver basis of the
Graver basis of $\Delta_2 \times \Delta_2$ which
attains this bound is gotten from the following
representation of the zero matrix:
\begin{eqnarray*}
& 3 \cdot 
\begin{pmatrix}
\phantom{-} 1 &  -1 & \phantom{-} 0 \, \\
\phantom{-} 0 & \phantom{-} 0 & \phantom{-} 0\, \\
 -1 & \phantom{-} 1 & \phantom{-} 0 \,
\end{pmatrix}
\,\, + \,\,
2 \cdot 
\begin{pmatrix}
 -1 &  \phantom{-} 1 & \phantom{-} 0 \, \\
\phantom{-}  0 &  -1 & \phantom{-} 1 \, \\
\phantom{-}  1 &  \phantom{-} 0 & -1 \, 
\end{pmatrix}
\,\, + \,\,
1 \cdot 
\begin{pmatrix}
\phantom{-}  1 & \phantom{-}  0 & -1 \, \\
 -1 & \phantom{-}  1 &\phantom{-}  0 \, \\
\phantom{-}  0 &  -1 & \phantom{-} 1\,  
\end{pmatrix}
\\ & \,\, + \,\,\,\,
2 \cdot 
\begin{pmatrix}
 -1 &  \phantom{-} 0 & \phantom{-} 1 \, \\
\phantom{-}  0 & \phantom{-}  1 & -1 \, \\
\phantom{-}  1 &  -1 &\phantom{-}  0 \, 
\end{pmatrix}
\,\,\,\, + \,\,\,\,
1 \cdot 
\begin{pmatrix}
\phantom{-}  0 & \phantom{-}  1 & -1 \, \\
\phantom{-}  1 &  -1 &\phantom{-}  0 \, \\
 -1 & \phantom{-}  0 &\phantom{-}  1 \, 
\end{pmatrix}.
\end{eqnarray*}

  \section{\bf Geometric Combinatorics}
  \label{sec:geometriccombinatorics}

A non-zero table $u$ in $\L(\A^{(r)})$ is a {\em circuit}
if the entries of $u$ are relatively prime and
the support of $u$ is minimal with respect to inclusion.
We define the {\em circuit complexity} $\, c(\A)\,$ as the maximum type
of any table that is a circuit of some
higher Lawrence lifting $\A^{(r)}$.
Since the set of circuits of $\A^{(r)}$ is a subset of the
Graver basis of $\A^{(r)}$, by
  \cite[Proposition 4.11]{Sturmfels-book}, we have
$$ c(\A) \quad \leq \quad g(\A) . $$
Example  \ref{exm:pqr} shows that there is no bound for
$g(\A)$ in terms of $n$ and $d$ alone.
On the other hand, such a bound does exist for
the circuit complexity $c(\A)$:

\begin{theorem}
\label{thm:circuits}
The circuit complexity of $\A$ is bounded above by $n-d+1$.
\end{theorem}

We shall derive this theorem from the following lemma, which can be  
rephrased as
``the circuits of any $\A^{(r)}$ are circuits among the circuits of  
$\A$''.

\begin{lemma}
\label{discirc}
Let $\C$ be the configuration consisting of all circuits of $\A$.
The non-zero rows of any circuit of type at least 3 of $\A^{(r)}$
are distinct (and not opposite) elements of $\C$ multiplied by
numbers which form a circuit of $\C$.
\end{lemma}

\begin{proof}
Let $u\in\L(\A^{(r)})$ represent a circuit of $\A^{(r)}$.
If two rows are opposite, then these two rows have the
sign pattern of an element of $\L(\A^{(2)})$,
and all other rows must be zero.
If some row $u^{(i)}$ is not a multiple of a circuit of $\A$,
then, by \cite[Lemma 4.10]{Sturmfels-book}, $u^{(i)}$ can be written
as  a non-negative rational conformal combination of circuits.
We can write
  $\,\alpha_0 u^{(i)} \,=\, \alpha_1 c_1 + \cdots + c_k  \alpha_k\,$
where the $\alpha_j$'s are positive integers and each $c_j\in\L(\A)$
is a circuit conformal to $u^{(i)}$. Then,
$\alpha_0 u$ can be decomposed as a sum of tables with support
strictly contained in that of $u$, a contradiction.

Let us now write  $u^{(i)}=\alpha_i c_i$, with $c_i\in\C$.
The vector of coefficients (the $\alpha_i$'s)
lies in $\L(\C)$. Again by \cite[Lemma 4.10]{Sturmfels-book}, if it is  
not a circuit of $\C$ then it
can be decomposed as a
non-negative rational conformal combination of circuits of $\C$.
As before, this decomposition
translates into a decomposition of some multiple of $u$ into tables
with strictly smaller support.
\end{proof}

\noindent {\sl Proof of Theorem \ref{thm:circuits}: }
The configuration $\C$ of Lemma \ref{discirc} has rank $n-d$.\qed


\medskip

Recall (e.g.~from \cite{OMbook} or \cite{Ziegler}) that
the {\em oriented matroid} of $\,\A^{(r)}\,$
is specified by the collection of all sign patterns
of circuits of $\A^{(r)}$.
Theorem \ref{thm:circuits} implies:

\begin{corollary}
If $c = c(\A) < r $ then the oriented matroid
of the higher Lawrence lifting  $\, \A^{(r)}\,$
is determined by the oriented matroid of $\, \A^{(c)}$.
\end{corollary}

The convex hull $\, {\rm conv}(\A^{(r)})\,$
of the higher Lawrence configuration
$\A^{(r)}$ is a {\em convex polytope}
of dimension $\,dr+n-d-1 \,$
in $\reals^{dr+n}$. A subset $C$ of $\A^{(r)}$
is a {\em face} of $\A^{(r)}$
if there exists a linear functional $\ell$
on $\reals^{dr+n}$ whose minimum over
$\,\A^{(r)} \,$ is attained precisely at the subset $C$.
Equivalently, the convex polytope $\, {\rm conv}(\A^{(r)})\,$
has a (geometric) face $F$ such that
$\, F \, \cap \, \A^{(r)}\,= \, C$.

\begin{corollary}
If $c = c(\A) < r $ then the set of faces of
$\A^{(r)}$ is determined by the set of faces of
$\,\A^{(c)}$.
In particular, a subset of $\A^{(r)}$ is a  face if and
only if its restriction to any
subtable with only $c$ rows is a  face of $\A^{(c)}$.
\end{corollary}

\begin{proof}
We use oriented matroid arguments as in \cite[\S 9]{OMbook}.
The faces of $\A^{(r)}$ are the complements
of the positive covectors of $\A^{(r)}$. Now an $r \times n$-table
of signs is a covector of $\A^{(r)}$ if and only if it is
orthogonal (in the combinatorial sense of \cite[\S 3]{OMbook})
  to all circuits of $\A^{(r)}$.
Since circuits have type at most $c$, the orthogonality relation
can be tested by restricting to subtables with $c$ rows only.
Hence an $r \times n$-table of signs is a (positive) covector
of $\A^{(r)}$ if and only if
every $c \times n$-subtable is  a (positive) covector
of $\A^{(c)}$.
\end{proof}

A basic result concerning the ``classical'' Lawrence
construction is that the oriented matroid of $\A$
can be recovered from the set of faces of $\A^{(2)}$,
and vice versa. This statement is
no longer true for higher Lawrence liftings.

\begin{example}
   \label{exm.mother}
We present a configuration $\A$ which has the property that
the set of faces of $\A^{(3)}$ cannot be recovered
from the oriented matroid of $\A$.
  Let  $d = 3$, $n = 6$ and consider the configurations
  \begin{equation}
  \A \, =\,
  \begin{pmatrix}
  4 & 0 & 0& 2& 1 & 1 \\
  0 & 4 & 0& 1& 2 & 1 \\
  0 & 0 & 4& 1& 1 & 2 \\
  \end{pmatrix} \quad \hbox{and}
\quad \A' \, =\,
  \begin{pmatrix}
  6 & 0 & 0& 2& 1 & 1 \\
  -1 & 4 & 0& 1& 2 & 1 \\
  -1 & 0 & 4& 1& 1 & 2 \\
\end{pmatrix}.
  \end{equation}
These two matrices determine the same
oriented matroid on $\{1,2,3,4,5,6\}$.
The configuration $\A^{(3)}$ consists of the $18$
column vectors of a $15 \times 18$-matrix formed as in (\ref{Athree}).
These vectors are indexed by the entries of a $3 \times 6$-table.
The following table is a positive covector of
the oriented matroid of $\A^{(3)}$:
\begin{equation}
\label{covector}
\begin{pmatrix}
+ & 0 & 0 & + & 0 & 0 \\
0 & + & 0 & 0 & + & 0 \\
0 & 0 & + & 0 & 0 & + \\
\end{pmatrix}
\end{equation}
To see this, multiply the $15 \times 18$-matrix $\A^{(3)}$
on the left by the vector
$$
\ell \quad = \quad
  \bigl(\, 3, -1, -1;\, -1, 3, -1;\, -1, -1, 3;\, 4,4,4,0,0,0 \,\bigr).
$$
This vector supports a face of $\,{\rm conv}(\A^{(3)})$,
and the elements of  $\A^{(3)}$ on that face
are indexed by the twelve zeros in (\ref{covector}).
We claim that the sign table (\ref{covector})
is not a covector of $(\A')^{(3)}$. If it were,
then there exists an analogous vector $\ell'$
such that $\ell' \cdot (\A')^{(3)}$ has
the same support as $\ell \cdot \A^{(3)}$.
This requirement leads to an inconsistent system
of linear equations for $\ell'$.
We conclude that while $\A$ and $\A'$
share the same rank $3$ oriented matroid, the
two polytopes $\, {\rm conv}(\A^{(3)}) \,$
and  $\, {\rm conv}( (\A')^{(3)}) \,$
are not combinatorially isomorphic.
  \end{example}

  \end{document}